\documentclass{amsart}
\newtheorem{Thm}{Theorem}[section]
\newtheorem{prop}[Thm]{Proposition}
\newtheorem {Lem}[Thm]{Lemma}
\newtheorem{Cor}[Thm]{Corollary}

\theoremstyle{remark}
\newtheorem{Rem}[Thm]{Remark}
\numberwithin{equation}{section}

\begin{document}

\title[First $L^p$-cohomology of groups]
{The first $L^p$-cohomology of some finitely generated groups and $p$-harmonic functions}

\author[M. J. Puls]{Michael J. Puls}
\address{Department of Mathematics \\
Eastern Oregon University \\
One University Boulevard \\
LaGrande, OR 97850 \\
USA}
\email{mpuls@eou.edu}

\begin{abstract}
Let $G$ be a finitely generated infinite group and let $p>1$. In this paper we make a connection between the first $L^p$-cohomology space of $G$ and $p$-harmonic functions on $G$. We also describe the elements in the first $L^p$-cohomology space of groups with polynomial growth, and we give an inclusion result for nonamenable groups.
\end{abstract}


\date{May 7, 2004}
\maketitle

\section{Introduction}
In this paper $G$ will always be a finitely generated infinite group and $S$ will always denote a symmetric generating set for $G$. Let $M$ be a right $G$-module and let $C^n(G,M)$ be the set of functions from $G^n = G\times \cdots \times G (n\mbox{ copies})$ to $M$. We now have a chain complex
$$0\longrightarrow M \stackrel{\delta_0}{\longrightarrow} C^1(G,M)\stackrel{\delta_1}{\longrightarrow} C^2(G,M)\stackrel{\delta_2}{\longrightarrow} \cdots \stackrel{\delta_{n-1}}{\longrightarrow} C^n(G,M)\stackrel{\delta_n}{\longrightarrow} \cdots $$
where for $n \geq 0$ 
$$(\delta_n f)(g_1,\ldots,g_{n+1}) = \left( f(g_2, \ldots, g_{n+1})\right)\cdot g_1 +$$
$$\sum_{k=1}^n (-1)^k f(g_1, \ldots, g_kg_{k+1}, \ldots, g_{n+1}) + (-1)^{n+1}f(g_1,\ldots ,g_n), $$
where $g_kg_{k+1}$ occupies the $k$th position of $f$. Each $\delta_n$ is a linear map and a calculation shows $\delta_n\delta_{n-1}=0$. The $n$th group cohomology of $G$ with coefficients in $M$ is denoted by $H^n(G,M)$ and is equal to $\mbox{ker } \delta_n/ \mbox{im } \delta_{n-1}$, where $\mbox{ker }\delta_n$ denotes the kernel of $\delta_n$ and $\mbox{im }\delta_{n-1}$ is the image of $\delta_{n-1}$ in $C^n(G,M)$. If $M$ is a Banach space then we can give $C^n(G,M)$ the compact open topology. This means $f_k\rightarrow f$ in $C^n(G,M)$ if and only if $f_k(g_1,\ldots ,g_n) \rightarrow f(g_1,\ldots,g_n)$. Note that $\mbox{ker } \delta_1$ is a Banach space. In general $\mbox{im } \delta_{n-1}$ is not closed in $\mbox{ker } \delta_n$. Let $\overline{\mbox{im }\delta_{n-1}}$ denote the closure of $\mbox{im } \delta_{n-1}$ in $\mbox{ker }\delta_n$. The quotient space $\overline{H}^n(G)=\mbox{ker } \delta_n/\overline{\mbox{im } \delta_{n-1}}$ is called the $n$th reduced cohomology space of $M$.

Let $\mathcal{F}(G)$ denote the set of real-valued functions on $G$. Let $f\in \mathcal{F}(G)$ and let $x\in G$, the right translation of $f$ by $x$ is the function defined by $f_x(g)=f(gx^{-1})$. For a real number $p\geq 1, L^p(G)$ will consist of those $f\in \mathcal{F}(G)$ for which $\sum_{g\in G} |f(g)|^p<\infty$, and $C_0(G)$ will consist of those $f\in \mathcal{F}(G)$ for which the set $\{ g\mid |f(g)| > \epsilon\}$ is finite for each $\epsilon > 0$. The sets $\mathcal{F}(G), C_0(G)$ and $L^p(G)$ are $G$-modules under right translations. In this paper we study $H^1(G,L^p(G))$, the first group cohomology of $G$ with coefficients in $L^p(G)$. We also study $\overline{H}^1_{(p)}(G)$, the first reduced $L^p$-cohomology space.  

This work was supported by a grant from the research office of Eastern Oregon University. The author would like to thank the university for their kind support.

\section{Preliminaries}\label{Preliminaries}
Let $\mathbb{R}G$ be the group ring of $G$ over $\mathbb{R}$. For $h \in \mathbb{R}G$ and $f \in \mathcal{F} \left( G \right)$ we define a multiplication from $\mathcal{F}\left( G \right) \times \mathbb{R}G \mbox{ into } \mathcal{F}\left( G \right)$ by
$$ (f \ast h)(g) =  \sum_{x\in G}f(gx^{-1})h(x).$$
Observe that $(f\ast (s-1))(g) = f(gs^{-1})- f(g)$ for each $g\in G$ and each $s\in S$. For $1\leq p \in \mathbb{R}$, let $D^p(G) = \{ f \in \mathcal{F}\left( G \right) \mid f \ast (s-1) \in L^p\left( G \right) \mbox{ for each } s \in S\}$. Recall that $S$ is a symmetric set of generators for $G$. We define a norm on $L^p\left( G \right)$ by $\parallel f \parallel_p= ( \sum_{g\in G}|f(g)|^p )^{\frac{1}{p}}$, where $f \in L^p\left( G \right)$. Let $h \in D^p\left( G \right)$ and let $e$ be the identity element of $G$. We define a norm on $D^p\left( G \right)$ by $\parallel h \parallel_{D^p(G)} = \left( \sum_{s\in S} \parallel h \ast (s-1)\parallel_p^p + |h(e)|^p \right)^{\frac{1}{p}}$. Under this norm $D^p\left( G \right)$ is a Banach space. Let $f_1$ and $f_2$ be elements of $D^p\left( G \right).$ We will write $f_1 \simeq f_2$ if $f_1 - f_2$ is a constant function. Clearly $\simeq$ is an equivalence relation on $D^p\left( G \right)$. Identify the constant functions on $G$ with $\mathbb{R}$. Now $D^p\left( G \right) / \mathbb{R}$ is a Banach space under the norm induced from $D^p\left( G \right)$. That is, if $[ f ]$ is an equivalence class from $D^p\left( G \right)/ \mathbb{R}$ then $\parallel [f] \parallel_{D^p\left( G \right)/ \mathbb{R}} = \left( \sum_{s \in S} \parallel f \ast (s-1) \parallel_p^p \right)^{\frac{1}{p}}$. We shall write $\parallel f \parallel _{D(p)}$ for $\parallel [f] \parallel_{D^p(G)/\mathbb{R}}$. The norm for $D^p(G)$ and $D^p(G)/\mathbb{R}$ depends on the symmetric generating set $S$, but the underlying topology does not. If $X \subset D^p(G)$, then $(\overline{X})_{D^p(G)}$ will denote the closure of $X$ in $D^p(G)$. Similarly if $Y \subset D^p(G)/\mathbb{R}$, then $(\overline{Y})_{D(p)}$ will denote the closure of $Y$ in $D^p(G)/\mathbb{R}.$ The cardinality of a set $A$ will be denoted by $|A|$.

For the chain complex
$$0\longrightarrow L^p(G) \stackrel{\delta_0}{\longrightarrow} C^1(G,L^p(G))\stackrel{\delta_1}{\longrightarrow} C^2(G,L^p(G))\stackrel{\delta_2}{\longrightarrow} \cdots $$
the map $\delta_0$ is given by $(\delta_0 f)(g)=f\ast(g-1)$. If $f\in \mathcal{F}(G)$ and $\delta_0 f \in C^1(G, L^p(G))$, then $\delta_1(\delta_0 f) =0$ which implies $D^p(G)/\mathbb{R} \subseteq \mbox{ ker }\delta_1$ since $f \ast (g-1) \in L^p(G)$ for all $g \in G$. We now show that the reverse inclusion is also true.
\begin{Lem}
Let $G$ be a finitely generated, infinite group. The kernel of $\delta_1:C^1(G,L^p(G)) \rightarrow C^2(G,L^p(G))$ is $D^p(G)/\mathbb{R}.$
\end{Lem}
\begin{proof}
Consider the chain complex
$$ 0\longrightarrow \mathcal{F}(G) \stackrel{\delta'_0}{\longrightarrow} C^1(G,\mathcal{F}(G)) \stackrel{\delta'_1}{\longrightarrow} C^2(G,\mathcal{F}(G)) \stackrel{\delta'_2}{\longrightarrow} \cdots$$
It was shown in \cite[Lemma 4.2]{houghton} that $H^1(G,\mathcal{F}(G))=0$. Thus for each $f\in \mbox{ker }\delta'_1$ there exists $\bar{f}$ in $\mathcal{F}(G)$ such that $\delta'_0\bar{f} = f.$

Let $f_1\in C^1(G,L^p(G)) \subseteq C^1(G,\mathcal{F}(G))$ and suppose $\delta_1f_1=0$. Then $\delta'_1f_1=0$ which implies there exists an $\bar{f}_1 \in \mathcal{F}(G)$ such that $\delta_0\bar{f}_1 = f_1$. The result now follows since $\bar{f}_1 \in D^p(G).$
\end{proof}

The map $\delta_0$ is an injection so we obtain the following:
\begin{enumerate}
\item[(a)] The first cohomology group of $G$ with coefficients in $L^p\left( G \right)$, denoted by $H^1\left( G, L^p\left( G \right)\right)$, is isomorphic with $D^p\left( G \right) / \left( L^p\left( G \right) \bigoplus \mathbb{R} \right).$
\item[(b)] The first reduced $L^p$-cohomology space of $G$, denoted by $\overline{H}_{(p)}^1 \left( G \right)$, is isometric with $D^p\left( G \right) / \overline{L^p\left( G \right) \bigoplus \mathbb{R}}$, where the closure is taken in $D^p\left( G \right).$
\end{enumerate} 

\section{Nonreduced $L^p$-cohomology and $p$-harmonic functions}
Let $f\in \mathcal{F}(G)$ and let $g\in G$. Suppose $1<p\in \mathbb{R}$ and define 
$$ \bigtriangleup_p f(g) := \sum_{s\in S} |f(gs^{-1})-f(g)|^{p-2}(f(gs^{-1})-f(g)).$$
In the case $1<p<2$, we make the convention that $|f(gs^{-1})-f(g)|^{p-2} (f(gs^{-1})-f(g))=0$ if $f(gs^{-1})=f(g).$ We shall say that $f$ is {\em $p$-harmonic} if $f\in D^p(G)$ and $\bigtriangleup_p f(g)=0$ for all $g\in G.$ Recall that $f$ is a harmonic function if $\sum_{s\in S} (f(gs^{-1})-f(g)) = 0$ for all $g \in G$. Let $HD^p(G)$ be the set of $p$-harmonic functions on $G$. Observe that the constant functions are in $HD^p(G)$. If $p=2$, then $HD^2(G)$ is the linear space of harmonic functions on $G$ with finite energy. In general, $HD^p(G)$ is not a linear space if $p\neq 2$. A wealth of information about $p$-harmonic functions on graphs and manifolds can be found in \cite{holopainen, holasoa2}. Many of the ideas in this section come from the paper \cite{yamasaki}. In this section we will give a decomposition of $D^p(G)$ that will allow us to make a connection between $p$-harmonic functions on $G$ and $\overline{H}^1_{(p)}(G)$. 

We begin by giving some preliminary results for both $D^p(G)$ and $HD^p(G)$. Let $f$ and $h$ be elements in $D^p(G)$ and let $1<p \in \mathbb{R}$. Define 
\begin{equation*}
\begin{split}
 \langle \bigtriangleup_ph, f\rangle  & := \sum_{g\in G}\sum_{s\in S} |h(gs^{-1})-h(g)|^{p-2} (h(gs^{-1})-h(g))(f(gs^{-1})-f(g))\\
 & = \sum_{g\in G} \sum_{s\in S} |(h\ast (s-1))(g)|^{p-2}((h\ast (s-1))(g))((f\ast (s-1))(g)).
\end{split}
\end{equation*}
The above sum exist since $\sum_{g\in G}\sum_{s\in S} \left| |h(gs^{-1})-h(g)|^{p-2} (h(gs^{-1}) - h(g))\right|^q < \infty$, where $\frac{1}{p} + \frac{1}{q}=1$. The next lemma will be used to help show the uniqueness of the decomposition of $D^p(G)$ that will be given in Theorem \ref{decompose}.
\begin{Lem} \label{helpuni}
Let $f_1$ and $f_2$ be functions in $D^p(G)$. Then $\langle \bigtriangleup_p f_1 - \bigtriangleup_p f_2, f_1-f_2 \rangle = 0$ if and only if $f_1\ast (s-1) = f_2\ast (s-1)$ for all $s\in S$.
\end{Lem}
\begin{proof}
Let $f_1, f_2 \in D^p(G)$ and assume there exists $s \in S$ such that $f_1\ast (s-1) \neq f_2\ast (s-1)$. Define a function $f:[0,1]\rightarrow \mathbb{R}$ by $f(t) = \sum_{g\in G}\sum_{s\in S} |f_1(gs^{-1}) - f_1(g) + t((f_2(gs^{-1}) - f_2(g)) - (f_1(gs^{-1}) - f_1(g)))|^p$. Observe that $f(0) = \parallel f_1 \parallel_{D(p)}$ and $f(1) = \parallel f_2 \parallel_{D(p)}$. A derivative calculation gives $\frac{df}{dt}\mid_{t=0} = p\langle \bigtriangleup_p f_1, f_2-f_1 \rangle = -p \langle \bigtriangleup_p f_1, f_1-f_2 \rangle$. It follows from Proposition 5.4 on page 24 of \cite{ekeland_temam} that $\parallel f_2 \parallel_{D(p)} > \parallel f_1 \parallel_{D(p)} - p \langle \bigtriangleup_p f_1, f_1-f_2 \rangle.$
Similarly, $\parallel f_1 \parallel_{D(p)} > \parallel f_2 \parallel_{D(p)} - p \langle \bigtriangleup_p f_2, f_2-f_1\rangle$. Hence, $p\langle \bigtriangleup_p f_1 - \bigtriangleup_p f_2, f_1-f_2 \rangle > 0$ if there exists $s \in S$ such that $f_1\ast (s-1) \neq f_2 \ast (s-1)$. 

Conversely, suppose $f_1\ast (s-1) = f_2 \ast (s-1)$ for all $s \in S$. Then $\langle \bigtriangleup_p f_1 - \bigtriangleup_p f_2, f_1-f_2 \rangle = 0$ since $f_1-f_2$ is a constant function on $G$.  
\end{proof}
For $g\in G$, define $\delta_g$ by $\delta_g(x)=0$ if $x\neq g$ and $\delta_g(g)=1$.
\begin{Lem} \label{vanish} Let $h\in \mathcal{F}(G)$. Then $h$ is a $p$-harmonic function if and only if $\langle \bigtriangleup_p h, \delta_x \rangle = 0$ for all $x\in G$.
\end{Lem}
\begin{proof}
Let $x\in G$ and suppose $h\in HD^p(G)$, then
$$ \langle \bigtriangleup_p h, \delta_x\rangle = -2 \sum_{s\in S} |h(xs^{-1}) - h(x)|^{p-2} (h(xs^{-1})-h(x))=0.$$
Conversely, if $\langle \bigtriangleup_p h, \delta_x\rangle =0$ for all $x\in G$, then $\sum_{s\in S} |h(xs^{-1})-h(x)|^{p-2} (h(xs^{-1})-h(x))=0$ for all $x\in G$.
\end{proof}
\begin{Rem} It follows immediately from the lemma that if $h\in HD^p(G)$, then $\langle \bigtriangleup_p h, f \rangle = 0$ for all $f\in \mathbb{R}G$.
\end{Rem}
\begin{prop} \label{orth}
If $h\in HD^p(G)$ and $f\in (\overline{L^p(G)})_{D^p(G)}$, then $\langle \bigtriangleup_p h, f \rangle = 0$. 
\end{prop}
\begin{proof}
Let $f\in (\overline{L^p(G)})_{D^p(G)}$ and let $h\in HD^p(G)$. There exists a sequence $\{ f_n \}$ in $\mathbb{R}G$ such that $\parallel f-f_n\parallel_{D^p(G)}\rightarrow 0$ as $n\rightarrow \infty$, since $(\overline{ \mathbb{R}G})_{D^p(G)} = (\overline{L^p(G)})_{D^p(G)}$. Now
\begin{equation*}
\begin{split}
0 & \leq |\langle \bigtriangleup_p h, f \rangle | \\
  & = |\langle \bigtriangleup_p h, f-f_n \rangle | \\
  & = | \sum_{g\in G} \sum_{s\in S} |(h\ast (s-1))(g)|^{p-2} ((h\ast (s-1))(g))(((f-f_n)\ast(s-1))(g))| \\
  & \leq \sum_{g\in G}\sum_{s\in S} | (h\ast (s-1))(g)|^{p-1} |((f-f_n)\ast (s-1))(g)| \\
  & \leq ( \sum_{g\in G}\sum_{s\in S} ( |(h\ast (s-1))(g)|^{p-1})^q ) \parallel f_n - f \parallel_{D(p)} \rightarrow 0
\end{split}
\end{equation*}
as $n\rightarrow \infty$. The last inequality follows from H\"older's inequality
\end{proof}

We now give a decomposition of $D^p(G)$ that will allow us to determine representatives for the nonzero classes in $\overline{H}^1_{(p)}(G)$.
\begin{Thm} \label{decompose}
Let $1< p \in \mathbb{R}$ and suppose $\overline{L^p(G)}_{D^p(G)} \neq D^p(G)$. If $f\in D^p(G)$, then we can write $f=u+h$, where $u \in (\overline{\mathbb{R}(G)})_{D^p(G)}$ and $h\in HD^p(G)$. This decomposition is unique up to a constant function.
\end{Thm}
\begin{proof}
Let $f \in D^p(G)$ and let $r$ equal the distance of $f$ from $(\overline{\mathbb{R}G})_{D^p(G)}$ in the $D^p(G)$-norm. Set $B=\{ v\in (\overline{\mathbb{R}G})_{D^p(G)} \mid \parallel f-v \parallel_{D^p(G)} \leq r+1 \}$. Now $B$ is a nonempty weakly compact set in the reflexive Banach space $D^p(G)$. The function $F(v)=\parallel f-v \parallel_{D^p(G)}$ is a weakly lower semi-continuous function on $B$, so $F(v)$ assumes a minimum value on $B$. This minimum must be $r$. Let $u \in (\overline{ \mathbb{R}G})_{D^p(G)}$ where $\parallel f-u \parallel_{D^p(G)} = r$ and set $h=f-u$. We now proceed to show that $h$ is $p$-harmonic on $G$. Let $t\in \mathbb{R}$ and let $w\in \mathbb{R}G$. Now, $\parallel h \parallel_{D^p(G)} = \parallel f-u \parallel_{D^p(G)} \leq \parallel f-(u - tw) \parallel_{D^p(G)}$ for all $t\in \mathbb{R}$. The minimum of $\parallel f -(u-tw) \parallel_{D^p(G)}$ occurs when $t=0$. Thus $\frac{d}{dt} ( \parallel h+tw \parallel_{D^p(G)})\mid_{t=0} = \sum_{g\in G}\sum_{s\in S} p|h(gs^{-1})-h(g)|^{p-2} (h(gs^{-1})- h(g))(w(gs^{-1})-w(g))=0$. Using $\delta_x$ for $w$ in the above derivative calculation we obtain $-2p\sum_{s\in S} |h(xs^{-1})-h(x)|^{p-2}(h(xs^{-1})-h(x))=0$. Thus $h$ is $p$-harmonic by Lemma \ref{vanish}.

We now show that this decomposition is unique up to a constant. Suppose $f= u_1+h_1 = u_2 + h_2,$ where $u_1,u_2 \in (\overline{\mathbb{R}G})_{D^p(G)}$ and $h_1, h_2 \in HD^p(G)$. Now, $\langle \bigtriangleup_p h_1 - \bigtriangleup_p h_2, h_1 - h_2 \rangle = \langle \bigtriangleup_p h_1 - \bigtriangleup_p h_2, u_2-u_1 \rangle = 0$ by Proposition \ref{orth} since $u_2- u_1 \in (\overline{\mathbb{R}G})_{D^p(G)}$. By Lemma \ref{helpuni}, $h_1\ast (s-1) = h_2\ast (s-1)$ for all $s\in S$. Thus $h_1-h_2$ is a constant function, which implies that $u_1-u_2$ is also a constant function.
\end{proof}
We saw in section \ref{Preliminaries} that $\overline{H}^1_{(p)}(G) = D^P(G) / \overline{ L^p(G) \bigoplus \mathbb{R}}$, so it follows from the theorem that each nonzero class in $\overline{ H}^1_{(p)}(G)$ can be represented by a nonconstant element in $HD^p(G)$. Thomas Schick has found an error in the proof of Theorem 5.3 from \cite{Puls}. We were unable to fix the proof of that theorem using the techniques of \cite{Puls}. The following corollary can be considered a partial fix of the theorem since many groups that have a central element of infinite order are also of polynomial growth.
\begin{Cor}
If $G$ is a finitely generated group with polynomial growth, then $\overline{H}^1_{(p)}(G) = 0$ for $1<p \in \mathbb{R}.$
\begin{proof}
It was shown in \cite[Corollary 1.10]{HoloSoa} that $HD^p(G)=\mathbb{R}$ if $G$ has polynomial growth. The result now follows.
\end{proof}
\end{Cor}

\section{Nonamenable groups and $L^p$-cohomology}
In this section we will give some results concerning the first $L^p$-cohomology space of nonamenable groups. 

Let $A$ be a subset of a group $G$ and define $\partial A := \{ x\in A \mid \mbox{ there exists } s\in S \mbox{ with } xs \notin A\}.$ We shall say that a group $G$ is {\em amenable} if it has an exhaustion $G_1 \subset G_2 \subset \cdots, \cup_{k=1}^{\infty} G_k = G$ by finite subsets such that $\lim_{k\rightarrow \infty} \frac{|\partial G_k|}{|G_k|} = 0$. A group that is not amenable is said to be {\em nonamenable}. Our first result shows how amenability affects the way $L^p(G)$ sits inside $D^p(G)/\mathbb{R}$.
\begin{Thm} \label{gui}
Let $G$ be a finitely generated infinite group and let $1\leq p \in \mathbb{R}$. Then $L^p(G)$ is closed in $D^p(G)/\mathbb{R}$ if and only if $G$ is nonamenable.
\end{Thm}
\begin{proof}
Assume that $G$ is nonamenable and let $f\in (\overline{L^p(G)})_{D(p)}=(\overline{\mathbb{R}G})_{D(p)}$. Now there exists a sequence $\{f_n\}$ in $\mathbb{R}G$ such that $f_n\rightarrow f$ in the Banach space $D^p(G)/\mathbb{R}$. Thus $\{f_n\}$ is a Cauchy sequence in $D^p(G)/\mathbb{R}$. It was shown in  \cite{Gerl} that there exist a constant $C_p$ such that $\parallel u \parallel_p \leq C_p\parallel u \parallel_{D(p)}$ for all $u\in \mathbb{R}G$ if and only if $G$ is a finitely generated nonamenable group. Thus $\{f_n\}$ is also a Cauchy sequence in $L^p(G)$. So there exists a $\bar{f}\in L^p(G)$ such that $\parallel \bar{f}-f_n \parallel_p \rightarrow 0$ as $n\rightarrow \infty$. Since $L^p$-convergence implies pointwise convergence we have that $\parallel (\bar{f} - f_n)\ast (s-1)\parallel_p \rightarrow 0$ as $n\rightarrow \infty$ for each $s\in S$. Hence $\parallel \bar{f}-f_n \parallel_{D(p)}\rightarrow 0$ as $n\rightarrow \infty$. Therefore, $\bar{f} = f$.

Conversely, suppose there exists an exhaustion $G_1\subseteq G_2 \subseteq \cdots , \cup_{k=1}^{\infty} G_k = G$ by finite subsets such that $\lim_{k\rightarrow \infty} \frac{|\partial G_k|}{|G_k|}=0$. Let $\chi_k$ denote the characteristic function on $G_k$. Now define a function $f_k$ on $G$ by $f_k := \frac{\chi_k}{\sqrt[p]{|G_k|}}$. Note that $\parallel f_k \parallel_p =1$ for all $k$. Computing the $D^p(G)/\mathbb{R}$-norm of $f_k$ we obtain $\parallel f_k \parallel_{D(p)}^p = \sum_{s\in S} \parallel f_k\ast (s-1)\parallel_p^p \leq 2\frac{|\partial G_k|}{|G_k|}$. Thus, $\lim_{k\rightarrow \infty} \parallel f_k \parallel_{D(p)} \rightarrow 0$ and $\parallel f_k \parallel_p =1$ for all $k$. Therefore, if $G$ is amenable, then $L^p(G)$ is not closed in $D^p(G)/\mathbb{R}$.
\end{proof}
The next result follows immediately from the theorem.
\begin{Cor}
Let $1\leq p \in \mathbb{R}$. If $G$ is a finitely generated nonamenable group, then $H^1( G, L^p(G)) = \overline{H}^1_{(p)}(G)$ and $H^1(G, L^p(G))$ is a Banach space.
\end{Cor}
We will now use the theorem to show $\overline{H}^1_{(p)}(G) \neq 0$ for groups with infinitely many ends.
\begin{Cor}
Let $1<p \in \mathbb{R}$. If $G$ is a finitely generated group with infinitely many ends, then $\overline{H}^1_{(p)} (G) \neq 0$.
\end{Cor}
\begin{proof}
We will prove the corollary by constructing a nonconstant harmonic function, say $h$, on $G$ that is an element of $D^p(G)$. It will then follow from the maximum principle for harmonic functions that $h\notin C_0(G)$. By the theorem $h \notin (\overline{L^p(G)})_{D^p(G)}$ since $L^p(G) \subset C_0(G)$. Theorem \ref{decompose} then shows $h=u+v$, where $u\in \overline{L^p(G)}_{D^p(G)}$ and $v$ is a nonconstant $p$-harmonic function, thus $h$ will represent a nonzero class in $\overline{H}^1_{(p)}(G).$ We now proceed to construct $h$ by using a technique that was used in the proof of Theorem 4.1 of \cite{soardiwoess}.

Recall that $S$ is a symmetric generating set for $G$. Define an element $P$ in $\mathbb{R}G$ by $P=\frac{1}{|S|} \sum_{s\in S}s$. Let $f\in L^p(G)$ and define a bounded linear operator on $L^p(G)$ by $P(f)= f\ast P$. Observe that $f$ is a harmonic function on $G$ if $f\ast (P-I) = 0$, where $I$ is the identity operator on $L^p(G)$. Let $k\in \mathbb{N}$ and denote by $P^k$ multiplication of $P$ with itself $k$ times. Since $G$ is nonamenable, $\parallel P \parallel < 1$ in the operator norm for $1< p \in \mathbb{R}$ \cite{Gerl}. Thus $F= - \sum_{k=0}^{\infty} P^k$ is a bounded operator on $L^p(G)$. Note $F$ is the inverse of $P-I$ in the space of bounded linear operators on $L^p(G)$. Let $X$ denote the Cayley graph of $G$ with respect to a the generating set $S$. Thus the vertices of $X$ are the elements of $G$, and $g_1,g_2 \in G$ are joined by an edge if and only if $g_1= g_2s^{\pm 1}$ for some generator $s$. Remove a finite number of vertices and edges of $X$ to obtain two disjoint, infinite, connected components of $X$. Let $X_1$ and $X_2$ denote the components. Define a function $h_1$ on $G$ by
$$ h_1(g) = \left\{ \begin{array}{cl}
     2   &  \mbox{if } $g$ \mbox{ is a vertex in } X_1 \\
     1   &  \mbox{if } $g$ \mbox{ is a vertex in } X_2 \\
     0   &  \mbox{otherwise}
\end{array} \right. .$$
The support of $h_1\ast (P-I)$ is contained in $\partial X_1 \bigcup \partial X_2$, so $h_1 \ast (P-I) \in \mathbb{R}G$. Thus $F(h_1\ast (I-P)) \in L^p(G)$. Set $h_2 = F(h_1\ast (I-P))$ and let $h= h_1-h_2$. Note $h \in D^p(G)$ since $h_1 \in D^p(G)$ and $h_2 \in L^p(G)$. Now, $h\ast (P-I) = h_1\ast (P-I) - h_1\ast (P-I) = 0$, so $h$ is a harmonic on $G$.
\end{proof}

Alain Valette has pointed out a simpler, but different, proof of the above corollary. We now proceed to give a quick sketch of his proof. By combining Lemma 2 with the remark after Proposition 1 of \cite{Bekkaval} we have that the dimension of $H^1(G, \mathbb{R}G)$ over $\mathbb{R}$ is the number of ends of $G$ minus one. The corollary now follows from Theorem \ref{gui}. Valette's proof also shows that the corollary is true for $p=1$.

The situation becomes unclear if $G$ is nonamenable with one end. For example, $H^1(SL_n (\mathbb{Z}), L^2(SL_n (\mathbb{Z}))) = 0$ for $n\geq 3$ since $SL_n(\mathbb{Z})$ has property $T$ when $n\geq 3$ \cite{guichardet}. On the other hand, if $G$ is a fundamental group of a closed Riemann surface of genus at least 2, then $H^1(G, L^2(G)) \neq 0$ \cite{cheegrum}. A good deal of information about $H^1(G, L^2(G))$ can be found in \cite{Bekkaval}.

\section{A Description of $H^1(G,L^p(G))$}
In this section we will describe the nonzero elements of $H^1(G,L^p(G))$ for groups that have polynomial growth of (precise) degree $n > p$. These results are a generalization of results from Section 6 of \cite{Puls}, where $H^1(G,L^2(G))$ was discussed.

Let $d>1$. We shall say that $G$ satisfies condition $S_d$ if there exists a constant $C > 0$ such that $\parallel f \parallel_{\frac{d}{d-1}} \leq C \parallel f \parallel_{D(1)}$ for all $f\in \mathbb{R}G$. If $f\in \mathcal{F}(G)$ and $t\geq 1$, then $f^t$ will denote the function $f^t(g) = (f(g))^t$. The following was proved in \cite{Puls} but we include it here for completeness.
\begin{Lem} \label{mvth}
Let $G$ be a finitely generated group and let $t$ be a real number greater than or equal to 2. If $f$ is a non-negative real-valued function in $\mathcal{F}(G)$, then 
$$ \parallel f^t \parallel_{D(1)} \leq 2t\sum_{g\in G} f^{t-1}(g)\left( \sum_{s\in S} |( f \ast (s-1))(g)| \right).$$
\end{Lem}
\begin{proof}
Let $g\in G$ and let $s\in S$. It follows from the Mean Value Theorem applied to $x^t$ that $(r^t - s^t) \leq t(r^{t-1} + s^{t-1})(r-s)$ where $r$ and $s$ are real numbers with $0 \leq s \leq r$. Thus $|f^t(gs^{-1}) - f^t(g)| \leq t (f^{t-1}(g) + f^{t-1}(gs^{-1}))|f(gs^{-1}) - f(g)|.$ We now obtain $\parallel f^t \parallel_{D(1)} = \sum_{g\in G} \sum_{s\in S} | ( f^t \ast (s-1)) (g)|  \leq t\sum_{g\in G} \sum_{s\in S} ( f^{t-1}(g) + f^{t-1}(gs^{-1})) |f(gs^{-1}) - f(g)| = 2t\sum_{g\in G} \sum_{s\in S} f^{t-1} (g) | f(gs^{-1}) - f(g)|$
\end{proof}
The next proposition is a generalization of Proposition 6.2 of \cite{Puls}. I would like to thank Thomas Schick for showing me this generalization.
\begin{prop} \label{est}
Let $2 \geq p \in \mathbb{R}$ and let $d > p$. If $G$ satisfies condition $S_d$, then there is a constant $C'>0$ such that $\parallel f \parallel_{\frac{pd}{d-p}} \leq C' \parallel f \parallel_{D(p)}$ for all $f\in \mathbb{R}G.$
\end{prop}
\begin{proof}
Set $t=\frac{pd-p}{d-p}$. By property $S_d$, Lemma \ref{mvth} and H\"older's inequality we have (assuming without loss of generality that $f$ is non-negative).
\begin{equation*}
\begin{split}
\parallel f^{\frac{pd-p}{d-p}} \parallel_{\frac{d}{d-1}} & \leq C \parallel f^{\frac{pd-p}{d-p}}\parallel_{D(1)} \\
 & \leq 2C \left(\frac{pd-p}{d-p}\right) \sum_{g\in G} f^{\frac{d(p-1)}{d-p}}(g) \left( \sum_{s\in S}|\bigl( f \ast (s-1)\bigr)(g)|\right) \\
 & =2C\left(\frac{pd-p}{d-p}\right)\sum_{g\in G} \sum_{s\in S}  f^{\frac{d(p-1)}{d-p}}(g) |f(gs^{-1}) - f(g) | \\
 & \leq 2C\left(\frac{pd-p}{d-p} \right) \parallel f^{\frac{d(p-1)}{d-p}} \parallel_{\frac{p}{p-1}} \parallel f \parallel_{D(p)}. 
\end{split}
\end{equation*}
Observe $\parallel f^{\frac{pd-p}{d-p}}\parallel_{\frac{d}{d-1}} = \parallel f^{\frac{pd}{d-p}}\parallel_1^{\frac{d-1}{d}}$ and $\parallel f^{\frac{d(p-1)}{d-p}} \parallel_{\frac{p}{p-1}} = \parallel f^{\frac{pd}{d-p}} \parallel_1^{\frac{p-1}{p}}$. Substituting we obtain $\parallel f^{\frac{pd}{d-p}} \parallel_1^{\frac{d-1}{d}} \leq C' \parallel f^{\frac{pd}{d-p}} \parallel_1^{\frac{p-1}{p}} \parallel f \parallel_{D(p)}$. Dividing this inequality by $\parallel f^{\frac{pd}{d-p}} \parallel_1^{\frac{p-1}{p}}$ and noting that $\parallel f^{\frac{pd}{d-p}}\parallel_1^{\frac{d-p}{pd}} = ( \parallel f \parallel_{\frac{pd}{d-p}}^{\frac{pd}{d-p}} )^\frac{d-p}{pd}$ will yield the claim in the proposition.
\end{proof}
We shall say that a group $G$ satisfies condition $(IS)_d$ if $|A|^{d-1} < C|\partial A|^d$ for all finite subsets $A$ of $G$ and a positive constant $C$. Varopoulous proves the following proposition on page 224 of \cite{varopoulous}, also see Chapter 1.4 of \cite{woess}.
\begin{prop}
A finitely generated group $G$ satisfies condition $(IS)_d$ for some $d\geq 1$ if and only if it satisfies condition $S_d$.
\end{prop}
Our next result will show that each nonzero class in $H^1(G,L^p(G))$, where $G$ is a group with polynomial growth of (precise) degree $d>2$, can be represented by a function in $L^{p'}(G)$ for some fixed real number $p' > p$.
\begin{Thm}
Let $G$ be a finitely generated group with polynomial growth of (precise) degree $d$. If $d>p \geq 2$, then each nonzero class in $H^1(G, L^p(G))$ can be represented by a function from $L^{\frac{pd}{d-p}}(G)$.
\end{Thm}
\begin{proof}
Varopoulous proves in the papers \cite{varo2, varo3} that $G$ has polynomial growth of (precise) degree $d$ if and only if $G$ satisfies condition $(IS)_d$. Thus $G$  also satisfies condition $S_d$.

Let $1_G$ denote the constant function one on $G$. If $1_G \in (\overline{ \mathbb{R}G})_{D^p(G)}$, then there exists a sequence ${f_n}$ in $\mathbb{R}G$ such that $ \parallel 1_G - f_n \parallel_{D^p(G)} \rightarrow 0$ but $\parallel f_n \parallel_{\frac{pd}{d-p}} \not\rightarrow 0$ contradicting Proposition \ref{est}. Hence $(\overline{\mathbb{R}G})_{D^p(G)} \neq D^p(G)$. By Theorem \ref{decompose} each $f\in D^p(G)$ can be represented uniquely by $u+h$, where $u \in (\overline{\mathbb{R}G})_{D^p(G)}$ and $h \in HD^p(G)$. By \cite[Corollary 1.10]{HoloSoa}, $HD^p(G) = \mathbb{R}$. Thus nonzero classes in $H^1(G, L^p(G))$ can be represented by functions in $(\overline{\mathbb{R}G})_{D^p(G)}\setminus L^p(G)$. Let $f \in (\overline{ \mathbb{R}G})_{D^p(G)}$, so there exists a sequence $\{ f_n \}$ in $\mathbb{R}G$ such that $f_n \rightarrow f$ in the Banach space $D^p(G)$. Observe that $f_n \rightarrow f$ in $D^r(G)$ for $p \leq r < \infty$. By Proposition \ref{est} $\{ f_n \}$ forms a Cauchy sequence in $L^{\frac{pd}{d-p}}(G)$. Let $\bar{f}$ be the limit of this sequence in   $L^{\frac{pd}{d-p}}(G)$. It now follows $f_n \rightarrow \bar{f}$ in $D^{\frac{pd}{d-p}} (G)$ since $\parallel (\bar{f} - f_n) \ast (s-1) \parallel_{\frac{pd}{d-p}} \rightarrow 0$ as $n \rightarrow \infty$ for each $s \in S$. Therefore $\bar{f} = f$ since $f_n \rightarrow f$ in $D^{\frac{pd}{d-p}}(G).$
\end{proof}
\begin{Cor}
Let $2<p \in \mathbb{R}$ and let $d$ be an integer greater than $p$. Each nonzero class in $H^1( \mathbb{Z}^d, L^p(\mathbb{Z}^d))$ can be represented by a function from $L^{\frac{pd}{d-p}}( \mathbb{Z}^d).$
\end{Cor}

\section{Some inclusion results}
Since $L^p(G) \subseteq L^{p'}(G)$ for $1 < p \leq p' \in \mathbb{R}$, a natural question to ask is how does $H^1(G,L^p(G))$ relate to $H^1(G, L^{p'}(G))$. In this section we will give some answers to this question.
\begin{Lem} \label{maxprin} Let $h \in HD^p(G)$. If $h \in C_0(G)$, then $h = 0$.
\end{Lem}
\begin{proof}
The set $\{ g \mid |h(g)| > \epsilon \}$ is finite for a given $\epsilon > 0$. Thus there exists an $x\in G$ such that  $|h(x) | \geq |h (g) |$ for all $g \in G$. It follows $h(x) = h (xs^{-1})$ for all $s \in S$ since $\sum_{s \in S} |h(xs^{-1}) - h(x) |^{p-2} h(xs^{-1}) = \sum_{g \in S} |h(xs^{-1})- h(x)|^{p-2} h(x)$. We now obtain $h(x) =  h (g)$ for all $g\in G$ since the Cayley graph of $G$ is connected. Therefore $h = 0$.
\end{proof}
We now give an inclusion result for nonamenable groups.
\begin{prop}
Let $G$ be a finitely generated nonamenable group. If $1 < p \leq p' \in \mathbb{R}$, then $H^1(G, L^p(G)) \subseteq H^1(G, L^{p'}(G))$.
\end{prop}
\begin{proof}
Let $f$ represent a nonzero class in $H^1(G, L^p(G))$. Thus $f \in D^p(G) \subseteq D^{p'}(G)$ and $f \notin L^p(G)$. By Theorem \ref{decompose} we can uniquely write $f = u +h$, where $u \in L^p(G)$ and $h$ is a nonconstant element in $HD^p(G)$. By Lemma \ref{maxprin} $h \notin C_0(G)$ so it follows $f \notin L^{p'}(G)$. Hence $f$ also represents a nonzero class in $H^1(G, L^{p'}(G))$. 
\end{proof}
We now finish this section by giving an example to show the above proposition is not true for amenable groups. Define $f \colon \mathbb{Z} \rightarrow \mathbb{R}$ by 
$$ f(n) = \left\{ \begin{array}{cl}
          \frac{1}{\sqrt[p]{n}} & n\geq 1 \\
            0                       & \mbox{otherwise}  \end{array} \right.  $$
for some $1 < p \in \mathbb{R}$. Observe that $f \notin L^p(\mathbb{Z})$ but $f \in L^{p'} (\mathbb{Z})$ for $p' > p.$ Now 
\begin{equation*}
\begin{split}
\sum_{n=1}^{\infty} |(f \ast (s-1)) (n) |^p & = \sum_{n=1}^{\infty} |f(n-1) - f(n)|^p \\ 
    & = 1+ \sum_{n=2}^{\infty} \left| \frac{\sqrt[p]{n} - \sqrt[p]{(n-1)}} {\sqrt[p]{n(n-1)}} \right|^p \\
            & \leq 1 + \sum_{n=2}^{\infty} \frac{1}{(n-1)^2}. 
\end{split}
\end{equation*}
Thus $f \in D^p(\mathbb{Z})$ which implies $f$ represents a nonzero class in $H^1( \mathbb{Z}, L^p(\mathbb{Z}))$ but $f$ is in the zero class of $H^1( \mathbb{Z}, L^{p'}( \mathbb{Z}))$ for $p' > p$.

\bibliographystyle{plain}
\bibliography{lpcopharm}
\end{document}